# DYNAMIC PROGRAMMING IN MIXED CONTINUOUS-DISCRETE TIME SCALES

S. A. Belbas
Mathematics Department
University of Alabama
Tuscaloosa, AL. 35487-0350. USA.

Abstract. We obtain the dynamic programming equations and optimality conditions akin to Pontryagin's extremum principle for certain mathematical models of hybrid control systems.

## 1. Introduction.

Many problems in control theory involve systems with both continuous-time and discrete-time components. This is typical for systems that admit a physical model in terms of a state that evolves in continuous time, and are also augmented by digital control devices that can be mathematically described by discrete-time equations. Consider, for example, a system of the following type:

$$\frac{dy(t)}{dt} = f(t, y(t), u(t), pz(t), qw(t));$$
$$z(\tau_{k+1}) = g(\tau_k, y(\tau_k), u(\tau_k), z(\tau_k), w(\tau_k))$$
(1.1)

Here, y and u are the continuous-time state and control, z and w are the discrete-time state and control, and p and q are interpolation operators that transform the discrete-time signals z and w into the continuous-time signals pz and qw. The discrete times $\{\tau_k\}$ may be fixed or they may be defined through analytical conditions involving t and x(t). One possible form of the operators p and q is

$$(pz(t))_i = \sum_{j=1}^{M_k} \varphi_{ijk}(t) z_j(\tau_k), \quad (qw(t))_i = \sum_{j=1}^{N_k} \psi_{ijk}(t) w_j(\tau_k)$$
$$\text{for } t \in (\tau_k, \tau_{k+1})$$
(1.2)

It will be useful to consider 2 versions of the continuous-time part and 2 versions of the discrete-time part; these versions are labelled (Ci) and (Dj), i,j=1,2:

$(C1). \quad \dfrac{dy(t)}{dt} = f(t, y(t), u(t), pz(t), qw(t))$

$(C2). \quad \dfrac{dy(t)}{dt} = f(t, y(t), u(t), pz(t))$

$(D1). \quad z(\tau_{k+1}) = g(\tau_k, y(\tau_k), u(\tau_k), z(\tau_k), w(\tau_k))$
$(D2). \quad z(\tau_{k+1}) = g(\tau_k, y(\tau_k), z(\tau_k), w(\tau_k))$

Every version of the type of hybrid systems we wish to consider can be specified as (Ci,Dj), for example (1.1) is (C1,D1).

It is well known (cf., e.g., Rozenvasser (1994)) that, for hybrid systems in general (not specifically the type of systems considered in this paper), it is necessary to take into account both, the discrete nature of part of the state and the control, and the continuous nature of the remaining part of the state and the control. It is generally insufficient to treat these types of systems either as purely discrete-time systems or as purely continuous-time systems.

These hybrid systems, as explained below, can be equivalently modelled as controlled ordinary differential equations with impulses, of the type treated in Bainov-Simeonov (1989). The associated problems of optimal control for such systems lead to impulsive first-order hyperbolic partial functional-differential equations; these equations are akin to (but technically different from) the equations satisfied by Lyapunov functions in Bainov-Simeonov (1989) and the impulsive partial differential equations studied in Bainov-Minchev (1995), Bainov-Kamont-Minchev (1995), and other works by the same group of authors.

The literature on control systems with impulses contains various models, problems, and methods, other than the ones presented in the present paper; we mention the following: the work of Zavalishchin and Sesekin (1991), which deals with generalized-function models and extremum conditions related to Pontryagin's principle; the work of Motta et al. (Motta-Rampazzo 1996, and further references therein), that relies on a method of imbedding an optimal control problem into a more general class of problems; the method of quasi-variational inequalities of Bensoussan and Lions (1982, 1987), originally for stochastic control problems and later extended to deterministic problems (cf. references in Bensoussan-Lions (1987)), in which method the times and magnitudes of impulsive control actions are determined on the basis of optimality criteria and the value function satisfies a single equation (or "quasi-variational inequality") or system, with nonlinearities of "max" or "min" type.

## 2. The basic case of dynamic programming for hybrid systems.

We consider a controlled system governed by a system of ordinary differential equations with impulses, as follows:

$$\dfrac{dx(t)}{dt} = f(t, x(t), u(t)) \quad for \ (t, x(t)) \notin \Gamma;$$
$$x(t^+) = x(t^-) + I(t, x(t), w(t)) \quad for \ (t, x(t^-)) \in \Gamma \qquad (2.1)$$

Here, $\Gamma$ is the impulse set in the (t,x)-space, and I is the impulse function. This model can incorporate certain control systems; for example, for the system (C2,D2), we can set

$$x(t) := [y(t)\ z(t)], \quad where$$
$$z(t) \equiv z(\tau_k),\ w(t) \equiv w(\tau\tau_k)\ for\ \tau_k < t < \tau_{k+1}; \tag{2.2}$$
$$I(t,x,w) := g(t,y,z,w) - z(t^-)$$

The set $\Gamma$ consists of the points $(\tau_k, x)$, $x \in I\!R$.

We consider the following optimal control problem: minimize a cost functional J given by

$$J := \int_0^T F(t, x(t), u(t))dt + \sum_{t:(t,x(t^-))\in\Gamma} \Phi(t, x(t^-), w(t)) + F_0(x(T^-)) \tag{2.3}$$

The parametrized functional is defined by

$$J_{s\xi}^{uw} := \int_s^T F(t, x_{s\xi}^{uw}(t), u(t))dt +$$
$$+ \sum_{t:t>s,(t,x_{s\xi}^{uw}(t^-))\in\Gamma} \Phi(t, x_{s\xi}^{uw}(t^-), w(t)) + F_0(x_{s\xi}^{uw}(T^-)) \tag{2.4}$$

where $x_{s\xi}^{uw}(t)$ solves

$$\frac{\partial x_{s\xi}^{uw}(t)}{\partial t} = f(t, x_{s\xi}^{uw}(t), u(t))\ for\ t > s\ and\ (t, x_{s\xi}^{uw}(t^-)) \notin \Gamma;$$
$$x_{s\xi}^{uw}(t^+) = x_{s\xi}^{uw}(t^-) + I(t, x_{s\xi}^{uw}(t^-), w(t))\ for\ t > s\ and\ (t, x_{s\xi}^{uw}(t^-)) \in \Gamma;$$
$$x_{s\xi}^{uw}(s^+) = \xi$$

$$(2.5)$$

The value function is defined by

$$V(s, \xi) := \inf_{u,w} J_{s\xi}^{uw}$$

$$(2.6)$$

and the corresponding equations of dynamic programming take the form

$$\inf_{a} \{\frac{\partial V(s,\xi)}{\partial s} + \sum_{i=1}^{n} \frac{\partial V(s,\xi)}{\partial \xi_i} f_i(s,\xi,a) + F(s,\xi,a)\} = 0 \quad for \ (s,\xi) \notin \Gamma;$$

$$V(s^-,\xi) = \inf_{b} \{V(s^+, \xi + I(s,\xi,b)) + \Phi(s,\xi,b)\};$$

$$V(T^-,\xi) = F_0(\xi)$$

(2.7)

## 3. Extensions of dynamic programming.

Let us consider a system of the following type:

$$\frac{dx(t)}{dt} = f(t,x(t),u(t)) \quad for \ (t,x(t^-)) \notin \Gamma;$$

$$x(t^+) = x(t^-) + I(t,x(t^-),u(t^-),w(t)) \quad for \ (t,x(t^-)) \in \Gamma$$

(3.1)

The corresponding functional to be minimized in this case is

$$J^{uw}_{as\ \xi} := \int_{s}^{T} F(t,x^{uw}_{as\ \xi}(t),u(t))dt + \sum_{t>s:(t,x^{uw}_{as\ \xi}(t^-))\in\Gamma} \Phi(t,x^{uw}_{as\ \xi}(t^-),u(t^-),w(t)) +$$

$$+ F_0(x^{uw}_{as\ \xi}(T))$$

(3.2)

where the additional subscript "a" refers to the requirement $u(t^-) = a$.
In this case the value function is defined as

$$V(s,\xi,a) := \inf_{u,w:u(t^-)=a} J^{uw}_{as\ \xi}$$

(3.3)

and the dynamic programming equations have the form

$$\inf_{\alpha} \{\frac{\partial V(s,\xi,a)}{\partial s} + \sum_{i=1}^{n} \frac{\partial V(s,\xi,a)}{\partial \xi_i} f_i(s,\xi,\alpha) + F(s,\xi,\alpha)\} = 0 \ ;$$

$$V(s^-,\xi,a) = \inf_{b} \{V(s^+, \xi + g(s,\xi,a,b)) + \Phi(s,\xi,a,b)\}; \ V(T^-,\xi,a) = F_0(\xi)$$

(3.4)

It is readily seen that the effect of the parameter "a" on the dynamic programming equations takes place exclusively through the impulse conditions, i.e. through the second equation in (3.4). This model includes the system (C1,D1).

Still another case may arise from aftereffect of the discrete control w on the impulsive part of the system. In that case, the impulsive condition becomes

$$x(t^+) = x(t^-) + I(t, x(t^-), w(t), w(\sigma(t)));$$

$$\sigma(t) := \sup\ \{\tau: \tau < t, (\tau, x(\tau^-)) \in \Gamma$$

(3.5)

In this case, we define

$$V(s, \xi, b) := \inf_{u, w\, :\, w(\sigma(s)) = b} J_{s\xi}^{uw}$$

(3.6)

Then the impulsive condition on V becomes

$$V(s^-, \xi, b) = \inf_c \{V(s^+, \xi + I(t, \xi, c, b), c) + \Phi(s, \xi, c, b)$$

(3.7)

These systems arise in problems of discrete-event control of continuous-time ordinary differential equations.

### 4. Connection with necessary conditions.

If we take the system of dynamic programming equations (2.7) and if we assume the existence of an optimal control function and associated trajectory, say $(u^*(t), w^*(t), x^*(t))$, then we must have

$$\frac{\partial V(s, x^*(s))}{\partial s} + \sum_{i=1}^{n} \frac{\partial V(s, x^*(s))}{\partial \xi_i} f_i(s, x^*(s), u^*(s)) + F(s, x^*(s), u^*(s)) = 0$$

for $(s, x^*(s^-)) \notin \Gamma$;

$$V(s^-, x^*(s^-)) = V(s^+, \xi + I(s, x^*(s^-), w^*(s))) + \Phi(s, x^*(s^-), w^*(s))$$

for $(s, x^*(s^-)) \in \Gamma$;

$$V(T^-, x^*(T^-)) = F_0(x^*(T^-))$$

(4.1)

Let us set

$$p_j^*(s) := \frac{\partial V(s, x^*(s))}{\partial \xi_j}, \quad p_j^*(s^\pm) := \frac{\partial V(s^\pm, x^*(s^\pm))}{\partial \xi_j};$$

$$H(s, x, p, u) := \langle p, f(s, x, u) \rangle + F(s, x, u);$$

$$K(s, x, p, w) := \langle p, I(s, x, w) \rangle + \Phi(s, x, w);$$

$$H^*(s, x, p) := H(s, x, p, u^*(s)); \quad K^*(s, x, p) := K(s, x, p, w^*(s))$$

(4.2)

Then differentiation of (4.1), under conditions of sufficient smoothness, leads to the Hamiltonian equations

$$\frac{dp_j^*(s)}{ds} = -\frac{\partial H^*(s, x^*(s), p^*(s))}{\partial \xi_j}, \quad \frac{dx_i^*(s)}{ds} = \frac{\partial H^*(s, x^*(s), p^*(s))}{\partial p_i}$$

for $(s, x^*(s)) \notin \Gamma$;

$$p_j^*(s^+) - p_j^*(s^-) = -\frac{\partial K^*(s, x^*(s^-), p^*(s^+))}{\partial \xi_j},$$

$$x_i^*(s^+) - x_i^*(s^-) = \frac{\partial K^*(s, x^*(s^-), p^*(s^+))}{\partial p_i}, \quad \text{for } (s, x^*(s^-)) \in \Gamma$$

(4.3)

and the dynamic programming equations (2.7) imply

$$H^*(s, x^*(s), p^*(s)) = \inf_u H(s, x^*(s), p^*(s), u) \text{ for } (s, x^*(s^-)) \notin \Gamma;$$

$$K^*(s, x^*(s^-), p^*(s^+)) = \inf_w K(s, x^*(s^-), p^*(s^+), w) \text{ for } (s, x^*(s^-)) \in \Gamma$$

(4.4)

The system (4.3, 4.4) is the impulsive analogue of the standard Hamiltonian equations and the corresponding extremum principle; H is the usual Hamiltonian, and K may be termed *the impulsive Hamiltonian.* We note that, in the case of a control acting linearly on the system, our Hamiltonian equations are consistent with the equations obtained in Zavalishchin-Sesekin (1991).

We can have an analogous set of equations for the problems of section 3. For example, for the first problem of section 3, we define

$$q^*_{aj}(s) := \frac{\partial V(s, x^*(s), a)}{\partial \xi_j}, \quad q^*_{aj}(s^\pm) := \frac{\partial V(s^\pm, x^*(s^\pm), a)}{\partial \xi_j};$$

$$H(s, x, p, u) := \langle p, f(s, x, u) \rangle + F(s, x, u);$$

$$K(s, x, p, a, w) := \langle p, I(s, x, a, w) \rangle + \Phi(s, x, a, w);$$

$$H^*(s, x, p) := \inf_u H(s, x, p, u); \quad K^*(s, x, p, a) := \inf_w K(s, x, p, a, w)$$

(4.5)

Then, with the insertion of the extra parameter "a", equations (4.3, 4.4) obtain.

## 5. Impulsive Riccati equations.

We consider a linear controlled system of the form

$$\frac{dx(t)}{dt} = P(t)x(t) + Q(t)u(t) \quad \text{for } (t, x(t^-)) \notin \Gamma;$$

$$x(t^+) = x(t^-) + M(t)x(t^-) + N(t)w(t) \quad \text{for } (t, x(t^-)) \in \Gamma$$

(5.1)

We shall treat a problem with quadratic costs as follows:

$$F(t, x, a) \equiv x^T A(t)x + 2x^T B(t)a + a^T C(t)a;$$

$$F_0(x) \equiv x^T A_0 x$$

$$\Phi(t, x, b) \equiv x^T \alpha(t)x + 2x^T \beta(t)b + b^T \gamma(t)b$$

(5.2)

We start with the Ansatz

$$V(s, x) = x^T K(s)x$$

(5.3)

We shall write $K^+$, $K^-$ for $K(s^+)$, $K(s^-)$, respectively. For $(t, x) \notin \Gamma$, we obtain a Riccati differential equation for K. For $(t, x) \in \Gamma$, we find the impulse conditions on K as follows: we calculate

$$V(s^+, x + I(s,x,b)) + \Phi(s,x,b) =$$
$$= x^T(E + M^T)K^+(E + M)x + 2x^T(E + M)K^+Nb + b^TN^TK^+Nb +$$
$$+ x^T\alpha x + 2x^T\beta b + b^T\gamma b$$

(5.4)

The necessary condition for minimizing the right-hand side of (5.4) gives

$$b^* = -(N^TK^+N + \gamma)^{-1} N^TK^+(E + M)x \equiv Lx$$

(5.5)

and the impulsive condition on the original dynamic programming system gives

$$K^- = K^+ + M^TK^+ + K^+M + M^TM + (E + M^T)K^+L + L^TK^+(E + M) + L^TN^TK^+NL$$
$$\text{for } (s,x) \in \Gamma$$

(5.6)

The impulse conditions (5.6) together with standard Riccati equations between impulses constitute the system of *impulsive Riccati equations.*